\documentclass[11pt]{article}
\usepackage{latexsym,amssymb,amsmath,amsfonts,amsthm}
\usepackage{graphics}
\usepackage{enumerate}
\usepackage{makecell}
\usepackage{graphicx}
\usepackage{algorithm}
\usepackage{algorithmic}
\renewcommand{\algorithmicrequire}{ \textbf{Input:}} 
\usepackage{mathrsfs}
\usepackage{subfigure}
\usepackage{color}
\usepackage{cite}
\usepackage{accents}
\usepackage{url}
\usepackage{multirow}

\makeatletter
\def\widebar{\accentset{{\cc@style\underline{\mskip14mu}}}}
\makeatother
\newtheorem{theorem}{Theorem}[section]

\newtheorem{assumption}[theorem]{Assumption}
\newtheorem{remark}[theorem]{Remark}
\newtheorem{definition}[theorem]{Definition}
\newtheorem{proposition}[theorem]{Proposition}

\definecolor{hw}{rgb}{0,0,0}
\begin{document}
\renewcommand{\theequation}{\arabic{section}.\arabic{equation}}
\title{\bf Deep unfolding as iterative regularization for imaging inverse problems }
\author{Zhuo-Xu Cui \thanks{Z.-X. Cui and Q. Zhu are contributed equally to this work.}
\thanks{Research Center for Medical AI, Shenzhen Institutes of Advanced Technology, Chinese Academy of Sciences, Shenzhen, China ({\tt zx.cui@siat.ac.cn}).}
\and
Qingyong Zhu\thanks{ Research Center for Medical AI, Shenzhen Institutes of Advanced Technology, Chinese Academy of Sciences, Shenzhen, China ({\tt qy.zhu@siat.ac.cn}).}
\and
Jing Cheng\thanks{Paul C. Lauterbur Research Center for Biomedical Imaging, Shenzhen Institutes of Advanced Technology, Chinese Academy of Sciences, Shenzhen, China ({\tt jing.cheng@siat.ac.cn}).}
\and
Dong Liang\thanks{Corresponding author. Research Center for Medical AI, Shenzhen Institutes of Advanced Technology, Chinese Academy of Sciences, Shenzhen, China and Paul C. Lauterbur Research Center for Biomedical Imaging, Shenzhen Institutes of Advanced Technology, Chinese Academy of Sciences, Shenzhen, China ({\tt dong.liang@siat.ac.cn}).}
}
\date{}

\maketitle

\begin{abstract}
Recently, deep unfolding methods that guide the design of deep neural networks (DNNs) through iterative algorithms have received increasing attention in the field of inverse problems. Unlike general end-to-end DNNs, unfolding methods have better interpretability and performance. However, to our knowledge, their accuracy and stability in solving inverse problems cannot be fully guaranteed. To bridge this gap, we modified the training procedure and proved that the unfolding method is an iterative regularization method. More precisely, we jointly learn a convex penalty function adversarially by an input-convex neural network (ICNN) to characterize the distance to a real data manifold and train a DNN unfolded from the proximal gradient descent algorithm with this learned penalty. Suppose the real data manifold intersects the inverse problem solutions with only the unique real solution. We prove that the unfolded DNN will converge to it stably. Furthermore, we demonstrate with an example of MRI reconstruction that the proposed method outperforms conventional unfolding methods and traditional regularization methods in terms of reconstruction quality, stability and convergence speed.

\vspace{.2in}


\end{abstract}

\section{Introduction}
\setcounter{equation}{0}

The inverse problem is commonly seen in the field of medical imaging. For example, the core goal of MRI is to reconstruct MR images from undersampled k-space data in order to shorten scan times and reduce artifacts from patient movement \cite{liang1992constrained}. Reduced exposure to excess x-ray radiation, the CT scan dose needs to be urgently reduced. Thus, the core goal of CT is to reconstruct CT images from undersampled angular projection data \cite{Buzug2011}.

The main difficulty in solving the inverse problem is that there are multiple (or infinite) solutions or the solutions are discontinuous to the measurement. Therefore, regularization is necessary for solving the inverse problem accurately and stably \cite{engl1996regularization}. Variational (Tikhonov-type) regularization \cite{Tikhonov1943OnTS} is one of the common methods. By solving a data consistency and a regularizer-composed variational problem, the inverse problem can be solved accurately and stably \cite{ito2014inverse}.
Early research mainly focused on designing handcraft regularizers to more accurately characterize the solution prior to the inverse problem. Successful handcraft regularizers include total variation (TV) \cite{rudin1992nonlinear}, wavelets \cite{mallat1999wavelet}, wavelet frames \cite{ron1997affine,dong2010mra}, etc. Since variational regularization is usually solved by iterative algorithms, iterative regularization has been developed as an attractive alternative regularization method. Iterative regularization usually consists of iterative algorithms with a termination criterion. The iterative termination yields stable approximations of the true solution when the measurement data contains noise \cite{KaltenbacherNeubauerScherzer+2008}. Classical successful iterative regularization methods include Landweber iteration \cite{hanke1995convergence}, Newton-type methods \cite{kaltenbacher1997some}, etc. However, the classical iterative regularisation method does not contain tools to characterize the solution prior, resulting in it typically performing poorly on imaging inverse problems.
To compensate for this deficiency, \cite{boct2012iterative,jin2014fast,jin2015inexact,jiao2016alternating,jiao2017preconditioned} modified the classical iterative scheme by drawing on the role of regularisers and introducing a penalty function to characterize the solution prior.

In numerous fields, deep learning (DL) has achieved great success in recent years. Given this, there is a focus on using DL to solve inverse problems in data-driven ways, such as DL-based CT, MRI, PET reconstruction \cite{wurfl2016deep,wang2016accelerating,xu2017200x,Su_2022,cui2022k,cui2022self,cao2022high}, etc. Early research mainly focused on using deep neural networks (DNNs) to directly learn the mapping between measurement to the solution of the inverse problem \cite{zhu2018image,zbontar2018fastmri}. Such methods lack the necessary interpretability and face the risk of instability \cite{antun2020instabilities}. As a result, the design of DL methods in the framework of regularisation has been considered. Some works have focused on using DNNs to learn functional regularisers, which are then brought into a variational regularization model to solve the inverse problem \cite{NEURIPS2018_d903e960,li2020nett,kobler2020total,mukherjee2020learned}. Particularly, Lunz et al. \cite{NEURIPS2018_d903e960} proposed a method for adversarial learning of regularisers to measure the distance to real data manifolds. Mukherjee et al. \cite{mukherjee2020learned} used an input-convex neural network (ICNN) \cite{amos2017input} to represent the adversarial regularizer. Because the adversarial regularizer meets convexity, the variational regularization model can be solved uniquely. However, because such methods rely heavily on traditional regularization models, in practice, many experiments have shown that it is difficult to outperform traditional regularization models by a significant margin.

On the other hand, a regularization model driven approach to network architecture design, namely the unfolding method, was proposed \cite{10.5555/3104322.3104374,NIPS2016_1679091c}. Specifically, an unfolded DNN is achieved by replacing certain modules in the iterative algorithm (e.g., proximal mapping) for variational regularization problems with learnable network modules. It is then trained in an end-to-end manner \cite{yang2018admm,adler2018learned,cheng2021learning,ke2021deep,ke2021learned,huang2021deep}. Since regularization models or iterative algorithms drive the network architecture, such methods are somewhat interpretable. In addition, numerical experiments have shown that such methods can achieve satisfactory performance for solving imaging inverse problems. However, the unfolding method breaks the architecture of the original iterative algorithm. Its output does not necessarily belong to the solution of the original regularization model \cite{cui2021equilibrated}. Then, their accuracy and stability in solving inverse problems cannot be fully guaranteed.

This paper attempts to bridge this theoretical gap for unfolding methods in solving inverse problems. To this end, we will draw on the techniques in traditional iterative regularization and functional regularizer learning methods to modify the training procedure of the unfolded DNNs. More precisely, we jointly learn a convex penalty function by an ICNN to characterize the distance to a real data manifold and train a DNN unfolded from the proximal gradient descent (PGD) algorithm with this learned penalty. For testing, we equip a new terminate criterion for the unfolded DNN and prove it theoretically as an iterative regularization method. Specifically, this work's main contributions and observations are summarized as follows.
\begin{itemize}
\item We propose a new procedure to jointly train a convex penalty function (to characterize the distance to a real data manifold) and train a DNN unfolded from a PGD algorithm with this learned penalty function simultaneously. Suppose the real data manifold intersects the inverse problem solutions with only the unique real solution. We prove that the unfolded DNN converges to the real solution when the measurement is clean.
\item When the measurement is noisy, we equip a new terminate criterion for the unfolded DNN and prove its finite-step-stopping property. Furthermore, when the measurement noise tends to zero, the output of the unfolded DNN is guaranteed to converge to the real solution of the inverse problem.
\item What's more, in the proposed training procedure, the loss function measures the distance from each iteration to the real data manifold rather than just the $L_2$-norm between the last iteration and the real data. Thus, the proposed method provides faster convergence numerically.
\item  Finally, we verify the proposed method on an example of MRI reconstruction. Experimental results show that our proposed method outperforms conventional unfolding and traditional regularization methods in terms of reconstruction quality and stability.
\end{itemize}

The remainder of the paper is organized as follows. Section \ref{sect_rw} reviews background. Section \ref{sect3} modifies the training procedure and proved that the unfolding method is an iterative regularization method. The implementation details are presented in Section \ref{sect4}. Experiments performed on several data sets are presented in Section \ref{sect5}. The last section \ref{sect6} gives some concluding remarks. All the proofs are presented in the Appendix.
	
\section{ Background}\label{sect_rw}
\subsection{Inverse Problems in Imaging}
The inverse problem considers the reconstruction of the unknown image $x$ from the measurement
\begin{equation}\label{eq:1}Ax=y\end{equation}
where $A:X\rightarrow Y$ is a bounded linear operator between two Hilbert spaces $X$ and $Y$, and $y\in Y$ is the measurement. The inner products and norms on $X$ and $Y$ will be simply denoted by $\langle\cdot,\cdot\rangle$ and $\|\cdot\|$, respectively, which should be clear from the context. For example, in MRI, $A$ is a subsampled Fourier transform composed of the Fourier transform and binary sampling operator \cite{liang1992constrained}. In CT, $A$ is a subsampled Radon transform which is a partial collection of line integrations \cite{Buzug2011}.
\subsection{Handcraft Regularization Methods}

Variational regularization is a successful and well-established method for solving inverse problems. Given a measurement $y$, a stable and exact solution to the inverse problem (\ref{eq:1}) is achieved by solving the following variational problem
 \begin{equation}\label{variational}\min_{x\in X}\frac{1}{2}\|Ax-y\|^2+\lambda r(x)\end{equation}
where $\|Ax-y\|^2$ is the data consistency, $r(x)$ is the regularizer which characterizes the solution prior and $\lambda$ is the regularization parameter whose choice crucially affects the performance of the method.
Since above model is usually solved by iterative algorithms, iterative regularization has been developed as an attractive alternative regularization method. Iterative regularization usually considers the inverse problem measurement containing noise, i.e., $y^\delta = y + n$, where $n$ denotes the noise with $\|n\|=\delta$. It consists of iterative algorithms with a terminate criterion, and the iterative termination yields stable approximations of the true solution. For example, classical Landweber iteration executes
 $$x^\delta_{k+1}=x^\delta_{k}-\eta A^*(Ax^\delta_{k}-y^\delta)$$
 When criterion $ \|Ax^\delta_{k}-y^\delta\|\leq\tau \delta$ ($\tau>0$) is satisfied, the above iteration terminates. However, the Landweber iteration does not contain tools to characterize the solution prior. It generally does not perform well in imaging inverse problems. To this end, we can solve the inverse problem by introducing a penalty function $f$ and using the proximal gradient descent method:
 \begin{equation}\label{ista}\begin{aligned}
\xi^\delta_{k+1}=&x^\delta_{k}-\eta A^*(Ax^\delta_{k}-y^\delta)\\
x^\delta_{k+1}\in&\arg\min_{x\in X}\frac{1}{2}\|x-\xi^\delta_{k+1}\|^2+f(x).
\end{aligned}\end{equation}
 It is worth mentioning that $f$ is a penalty function used to characterize the prior information of the solution $x$ and does not act as a regularizer.
\subsection{Unfolding Methods}
In the above methods, regularisers or penalty functions are handcrafted. It is not easy to accurately and adaptively characterize data priors. The corresponding algorithm unfolds the DNN by replacing the operators on regularisers or penalty functions with network modules. The unfolded DNN is then trained end-to-end to achieve data-driven prior learning. Specifically, for (\ref{ista}), we can replace the proximal operator about $f$ with the network module $\mathcal{S}_{\theta_{k}}$ with parameter $\theta_{k}$ at iteration $k$. Given training data $\{(x^n,y^n)\}_{n=1}^N$, the unfolded DNN reads:
 \begin{equation}\label{u_dnn}\begin{aligned}
\xi_{k+1}^n=&x_{k}^n-\eta A^*(Ax_{k}^n-y^n)\\
x_{k+1}^n=&\mathcal{S}_{\theta_{k+1}}(\xi_{k+1}^n).
\end{aligned}\end{equation}
When the network depth $K$ is determined, the network output $x_K^{n}$ can be considered a vector function of the parameter $\Theta:=[\theta_1,\ldots,\theta_K]$. Therefore, the DNN is trained by minimizing the distance between $x_K^{n}(\Theta)$ and the truth $x^n$.
\begin{equation}\label{loss2}\min_{\Theta}\frac{1}{N}\sum_{n=1}^N\|x_{K}^n(\Theta)-x^n\|^2.\end{equation}

It is easy to see that the unfolded DNN (\ref{u_dnn}) is closer in the form of the iterative regularization method (\ref{ista}). In this paper, we will modify the training procedure of the unfolded DNN (\ref{u_dnn}) and prove it is an iterative regularization.

\section{Unfolded DNN for Inverse Problems}\label{sect3}
In this section, we will redesign the training procedure of unfolded DNN (\ref{u_dnn}) and show that it will conform to the iterative regularization framework. Specifically, we will jointly train a convex penalty function $f(\cdot)$ and the unfolded DNN (\ref{u_dnn}). The training is done so that the convex penalty function measures the distance between each iteration $x_{k}^n$ and a real data manifold, and the module $\mathcal{S}_{\theta_{k}}$ in the unfolded DNN approximates the proximal mapping of the learned convex penalty function.
\subsection{Training Procedure}
Now we describe the jointly training procedure of the penalty function and the unfolded DNN in detail. Assumes that the training data $\{x^n\}_{n=1}^N\in X$ independent samples from distribution of real data (truth images) $\mathbb{P}_x$, $\{y^n\}_{n=1}^N\in Y$ independent samples from distribution of measurements $\mathbb{P}_y$. Based on the measurement $\{y^n\}_{n=1}^N$, the initial values $\{x_0^n\}_{n=1}^N$ of the unfolded DNN can be obtained by a pseudo-inverse mapping $A^{\dag}$ of $A$, i.e., $\{x_0^n\}_{n=1}^N=\{A^{\dag}y^n\}_{n=1}^N\sim\mathbb{P}_0=A^{\dag}_{\#}\mathbb{P}_y$. Based on these training data, the joint training procedure can be implemented by alternately optimizing the following two subproblems:

\textbf{Subproblem 1:} At the $t$-th loop, fixing the penalty function $f_{\phi^t}$, for each iteration $x^n_k$, we force the network module $\mathcal{S}_{\theta_{k}}$ to approximate the proximity mapping of $f_{\phi^t}$ at $\xi_{k}^n:=x^n_k-\eta A^*(Ax^n_k-y^n)$, $k=0,\ldots,K-1$. Thus, we defined the loss function $J_{1,t}([\theta_1,\ldots,\theta_K])$ (\ref{loss}) and minimized it using the Adam algorithm \cite{kingma2014adam}. It will be executed by $T_{\Theta}$ iterations, with step size $\gamma$ and decay rates $\beta_1,\beta_2$. Suppose $f_{\phi_t}$ is a convex function. In that case, the optimal network modules $\mathcal{S}_{\theta_{k}}$ can make (\ref{loss}) obtain a unique minimal value, which means that optimal $\mathcal{S}_{\theta_{k}}$ is the proximal mapping of $f_{\phi_t}$ at $\xi_k^n$.

\textbf{Subproblem 2:} Fixing the already updated $\{\mathcal{S}_{\theta_k^{t+1}}\}_{k=1}^K$, we minimize the loss function $J_{2,t}(\phi)$ (\ref{penalty}) to learn the penalty function $f_{\phi}$ adversarially. In particular, we use ICNN to represent $f_{\phi}$ to ensure its convexity, whose architecture will be shown in the following section. If the minimum value of the loss function (\ref{penalty}) is taken, we will show that $f_{\phi}$ measures the distance from iterates $\mathcal{S}_{\theta_k^{t+1}}(\xi^n_k)$ to the real data manifold.

It is worth mentioning that the conventional loss function (\ref{loss2}) only forces that the last iteration converges to the true solution. Given that $f_{\phi_t}$ represents the distance to the real data manifold, our proposed loss function (\ref{loss}) will force each iteration to converge to the true solution. Therefore, this training strategy can potentially make the unfolded DNN converge faster numerically.
\begin{figure*}[thbp]
\begin{center}
\subfigure{\includegraphics[width=0.95\textwidth,height=0.22\textwidth]{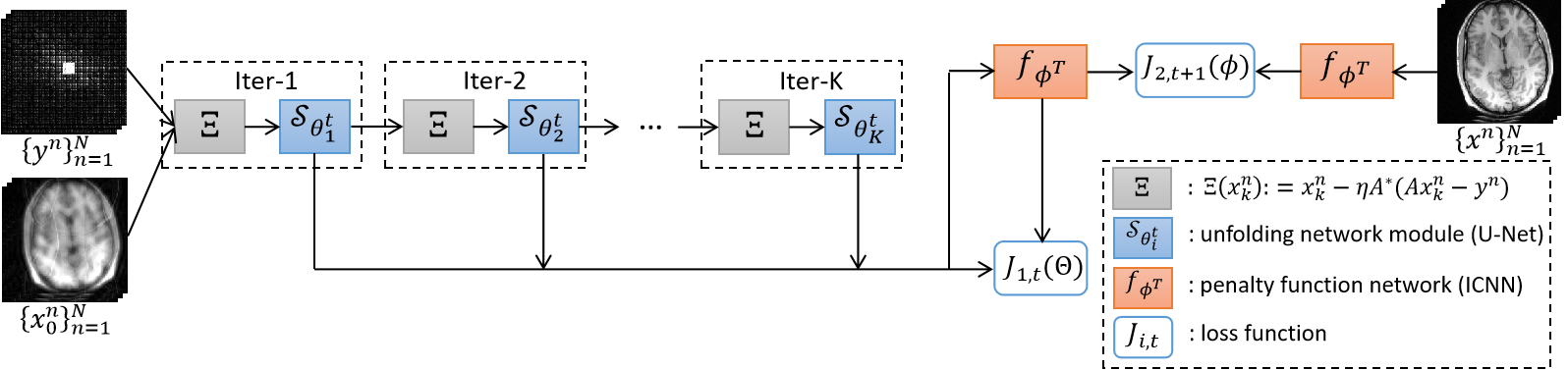}}
\end{center}
\caption{Illustration of the jointly training procedure of the unfolded DNN and penalty function.}
\label{f1}
\end{figure*}

The specific training steps are detailed in Algorithm \ref{alg:1}:
\begin{algorithm}[H]
	\caption{Training Procedure.}
	\label{alg:1}
	\begin{algorithmic}[1]
		\REQUIRE input $\{x^n\}_{n=1}^N\sim \mathbb{P}_x$, $\{y^n\}_{n=1}^N\sim \mathbb{P}_y$, Adam optimizer, positive constants $\eta$, $\gamma$, $\beta_1$, $\beta_2$ $K$, $T$, $T_{\Theta}$ and $T_{\phi}$.\\
        \renewcommand{\algorithmicrequire}{ \textbf{Initialize:}}
        \REQUIRE $\{x_0^n\}_{n=1}^N$, unfolded DNN $\mathcal{G}_{\Theta^0}: x_{k+1}^n=\mathcal{S}_{\theta_k^0}(\xi_{k}^n)$, $\xi_{k}^n:=x^n_k-\eta A^*(Ax^n_k-y^n)$, $k=0,\ldots,K-1$ with initial parameter $\Theta^0=[\theta_1^0,\ldots,\theta^0_K]$, convex and 1-Lipschitz continuous non-negative network $f_{\phi^0}(\cdot)$ with initial parameter $\phi^0$.\\
        \FOR{$t=0,\ldots,T-1$}
        \STATE {\bfseries Define loss function for $\mathcal{G}_{\Theta}$:} \begin{equation}\label{loss}J_{1,t}(\Theta)=\sum_{n=1}^N\sum_{k=1}^K\left[\frac{1}{2}\|\mathcal{S}_{\theta_k}(\xi_{k}^n)-\xi_{k}^n\|^2+f_{\phi^t}(\mathcal{S}_{\theta_k}(\xi_{k}^n))\right],~\xi_{k}^n=x_k^n-\eta A^*(Ax_k^n-y^n).\end{equation}
        \STATE Set $\Theta^t_0=\Theta^t$.
        \FOR{$s=0,\ldots,T_\Theta-1$}
        \STATE $$\Theta^t_{s+1}=\text{Adam}(\Theta^t_s,J_{1,t}(\Theta^t_s),\gamma,\beta_1,\beta_2).$$
        \ENDFOR
        \STATE $\Theta^{t+1}=\Theta^t_{T_\Theta}$.
        \STATE {\bfseries Define loss function for $f_{\phi}$:} \begin{equation}\label{penalty}J_{2,t+1}(\phi)=\mathbb{E}_{x^n\sim \mathbb{P}_x}[f_{\phi}(x^n)]-\frac{1}{K}\sum_{k=1}^K\mathbb{E}_{x^n_0\sim \mathbb{P}_0}\left[f_\phi(\mathcal{S}_{\theta_k^{t+1}}(\xi^n_k))\right].\end{equation}
        \STATE Set $\phi^t_{0}=\phi^{t}$.
        \FOR{$s=0,\ldots,T_\phi-1$}
        \STATE $$\phi^t_{s+1}=\text{Adam}(\phi^t_s, J_{2,t+1}(\phi^t_s),\gamma,\beta_1,\beta_2).$$
        \ENDFOR
        \STATE $\phi^{t+1}=\phi^t_{T_\phi}$.
        \ENDFOR
		\ENSURE $\mathcal{G}_{\Theta^T}$ (or $\{\mathcal{S}_{\theta_k^T}\}_{k=1}^K$) and $f_{\phi^T}$.
	\end{algorithmic}
\end{algorithm}

Now, let us normalize the assumptions and state the result on $f_{\phi}$. For Algorithm \ref{alg:1}, we can define a sequence of algorithmic operators $\mathcal{T}_k:X\rightarrow X$ to represent the iteration of unfolded DNN (\ref{u_dnn}). In particular, we define
$$ \mathcal{T}_{k+1}(u)= \mathcal{S}_{\theta_{k+1}}(\mathcal{T}_k(u)-\eta A^*(A\mathcal{T}_k(u)-y)).$$
Then the $k$th iteration can be reformulated as
$$x_{k}= \mathcal{T}_{k}(x_{0})$$
Initialized with distribution $\mathbb{P}_0$, we also define the distribution $\mathbb{P}_k$ as the push forward of the algorithmic operator $\mathcal{T}_{k}$ applied on $\mathbb{P}_0$, i.e.,
$$\mathbb{P}_k=(\mathcal{T}_{k})_{\#}\mathbb{P}_0.$$
The next assumption will relate $\mathbb{P}_k$ to the real data distribution $\mathbb{P}_x$. Before giving this assumption, we make some definitions and assumptions about the true data manifold $\mathcal{M}$.
Drawing on literatures \cite{NEURIPS2018_d903e960,doi:10.1137/20M1376790}, we first formalize the idea that real data are contained in a low dimensional manifold $\mathcal{M}$.
\begin{assumption}\label{assump3}
Assume the distribution $\mathbb{P}_x$ is supported on the convex compact set $\mathcal{M}$ of $X$, i.e., $\mathbb{P}_x(\mathcal{M}^c)=0$.
\end{assumption}
Based on the above assumption, the projection onto data manifold $\mathcal{M}$ can be defined as
$$P_{\mathcal{M}}:X\rightarrow \mathcal{M},~P_{\mathcal{M}}(x):= \arg\min_{z\in \mathcal{M}}\|z-x\|$$
This projection can be used to express the pointwise distance function
$$d_{\mathcal{M}}(x):=\inf_{z\in \mathcal{M}}\|z-x\|=\|P_{\mathcal{M}}(x)-x\|.$$
The strong assumption on each real data can be reconstructed from its measurement data is discarded. This paper follows \cite{doi:10.1137/20M1376790} and makes the following weaker assumption:
\begin{assumption}\label{assump4}
For any $k\leq K$, the push forward of the projection operator $P_{\mathcal{M}}$ on distribution $\mathbb{P}_k$ recovers the real data distribution $\mathbb{P}_x$ up to a set of measure zero, i.e., $\mathbb{P}_k$ pushed forward by $P_{\mathcal{M}}$ satisfies
$$(P_{\mathcal{M}})_{\#}\mathbb{P}_k=\mathbb{P}_x.$$
\end{assumption}
Drawing on the proof technique from \cite{NEURIPS2018_d903e960}, we show that if the loss function (\ref{penalty}) is obtained minimally, the trained convex penalty function can measures the distance from the iterates to the real data manifold.
\begin{theorem}\label{thm:1}Suppose Assumptions \ref{assump3} and \ref{assump4} hold. For any $k\leq K$, the minimizer to the functional
\begin{equation}\label{adv}\mathbb{E}_{X\sim \mathbb{P}_x}\left[f(X)\right]-\frac{1}{K}\sum_{k=1}^K\mathbb{E}_{Z^k\sim \mathbb{P}_k}\left[f(Z^k)\right]\end{equation}
is given by the distance of $Z^k\sim\mathbb{P}_k$ to the real data manifold $\mathcal{M}$, where $f$ is constrained to a convex and 1-Lipschitz continuous non-negative function.
\end{theorem}
In order to keep the gap between Algorithm \ref{alg:1} and Theorem \ref{thm:1} as small as possible, we will use ICNN to guarantee the convexity of penalty function $f_\phi$ and spectral normalization \cite{miyato2018spectral} or gradient penalty \cite{NIPS2017_892c3b1c} to characterize the continuity.

According to Theorem \ref{thm:1}, it is not difficult to understand how the proposed method works, i.e., the learned penalty function measures the distance onto a real data manifold, and network modules approximate the proximity mapping of this learned penalty. Therefore, the proposed method is an interpretable DL method.

\subsection{Testing Procedure}
Based on the learned proximal mapping $\{\mathcal{S}_{\theta_k^T}\}_{k=1}^K$ (the output of Algorithm \ref{alg:1}), we can execute the corresponding PGD algorithm to solve the inverse problem (\ref{eq:1}). However, in applications, measurements often contain noise. If the PGD algorithm is executed directly, the algorithm will show semi-convergence, i.e., it will start the iteration by traveling towards the true solution but will deviate from it as the iteration progresses \cite{jiao2016alternating,jiao2017preconditioned}. Therefore, it is important to terminate properly.

In this section, we consider the measurement containing $\delta$-level noise $n$, i.e., $y^\delta = y + n$ and $\|n\|=\delta$. Then, we equip a terminate criterion (\ref{terminate}), modified from the Morozov's criterion \cite{nair2003morozov}, to the learned PGD algorithm. It is worth noting that $f_{\phi^T}^*$ in the termination criterion may cause difficulties in practice. Theorem \ref{thm:1} shows that if $\phi^{T}$ obtains the minimum of (\ref{penalty}), we have $f_{\phi^T}^*=d_{\mathcal{M}}(x^\dag)=0$ where $x^\dag$ is the real image. If $\phi^{T}$ fails to obtain the minimum of (\ref{penalty}), $f_{\phi^T}(x^\dag)(\geq f_{\phi^T}^*)$ will also tend to be small value in the network training. Therefore, in practice, we can assume that $f_{\phi^T}^*\approx0$. Next, we will show that the terminated output solution satisfies regularity.
\begin{algorithm}[H]
	\caption{Testing Procedure.}
	\label{alg:2}
	\begin{algorithmic}[1]
		\REQUIRE learned network modules $\{\mathcal{S}_{\theta_k^T}\}_{k=1}^K$, penalty function $f_{\phi^T}$, $f^*_{\phi^T}=\min_{x\in \{x\in X|Ax=y\}}f_{\phi^T}(x)$, $y^{\delta}=y+n$ and positive constants $\tau$, $\eta$ and $K$.\\
		\renewcommand{\algorithmicrequire}{ \textbf{Initialize:}}
        \REQUIRE $x_0^\delta$.\\
        \FOR{$k=0,\ldots,K-1$}
        \STATE
        $$x_{k+1}^\delta=\mathcal{S}_{\theta_k^T}(x_k^\delta-\eta A^*(Ax_k^\delta-y^\delta)).$$
        \IF{\begin{equation}\label{terminate}\|Ax_{k+1}^\delta-y^\delta\|^2+f_{\phi^T}(x_{k+1}^\delta)-f_{\phi^T}^*\leq\tau^2\delta^2\end{equation}}
        \STATE break;
        \ENDIF
        \ENDFOR
		\ENSURE output $x_{k+1}^\delta$.
	\end{algorithmic}
\end{algorithm}

\subsection{Convergence and Regularity}
Before analyzing the convergence and regularity, we normalize the assumption on the learned network modules.
\begin{definition}We say $x$ solves minimization problem $\min_{z\in X}g(z)$ inexactly with error $\epsilon$ if $x$ satisfies $$g(x)\leq\min_{z\in X}g(z)+\epsilon$$\end{definition}
Algorithm \ref{alg:1} has difficulty ensuring that the loss function \ref{alg:1} converges to zero absolutely, resulting in difficulties for $\mathcal{S}_{\theta_k^T}(\xi_k)$ to accurately approximate the proximal mapping of $f_{\phi^{T}}$. Thus, we make the following more relaxed assumption.
\begin{assumption}\label{assump1}
In Algorithm \ref{alg:2}, the learned network modules $\mathcal{S}_{\theta_k^T}(\xi_k)$ solves $\min_z\frac{1}{2}\|z-\xi_k\|^2+f_{\phi^{T}}(z)$ inexactly with error $\epsilon_k$ and $\sum_{k=0}^{+\infty}\epsilon_k<+\infty$, where $\xi_k=x_k^\delta-\eta A^*(Ax_k^\delta-y^\delta)$.
\end{assumption}

\begin{remark} In practice, if Assumption \ref{assump1} does not hold, since $f_{\phi^{T}}$  is convex, we can obtain a solution satisfying Assumption \ref{assump1} by performing additional few (sub)gradient descent steps with $\mathcal{S}_{\theta_k^T}(\xi_k)$ as the initial value.
\end{remark}
Consider first that the measurement is mixed with noise. In this case, we will show that Algorithm \ref{alg:2} will be terminated after a finite number of iteration steps.

\begin{proposition}\label{pro:1}Suppose Assumption \ref{assump1} hold. Let $0<\eta<\frac{1}{2}$ and $\tau>\max\{\frac{3}{2-\eta},\frac{\eta}{2}\}$. When the measurement is mixed with noise, i.e., $\delta>0$, Algorithm \ref{alg:2} will terminate after finite $k_*$ steps of iteration.\end{proposition}

When the measurements are clean, i.e., $\delta = 0$, we will normalize some assumptions and state the result on the convergence of Algorithm \ref{alg:2}.
\begin{definition}
An solution $x^*\in X$ is called a $f$-minimizing solution to inverse problem (\ref{eq:1}) if it satisfies
$$f(x^*)\leq f(x),~\forall x\in\{u\in X|Au=y\}.$$
\end{definition}

\begin{assumption}\label{assump2}
The real data manifold $\mathcal{M}$ intersects the solutions of inverse problem (\ref{eq:1}) with the unique real solution, i.e., $\mathcal{M}\cap\{x\in X|Ax=y\}=\{x^{\dag}\}$
\end{assumption}
\begin{remark}If $f$ measures the distance to the real data manifold $\mathcal{M}$, i.e., $f(\cdot)=d_{\mathcal{M}}(\cdot)$, and Assumption \ref{assump2} holds, then the $f$-minimizing solution is the unique real solution to the inverse problem (\ref{eq:1}).\end{remark}

Assumption \ref{assump2} is strong, and we will discuss the convergence of Algorithm \ref{alg:2} when it holds or does not hold, respectively.
\begin{theorem}\label{thm:2}~\\
(a) Consider $\phi^{T}$ obtaining the minimum value of (\ref{penalty}). Suppose Assumption \ref{assump3}, \ref{assump4}, \ref{assump1} and \ref{assump2} hold, and the measurement is clean. The sequence generated by Algorithm \ref{alg:2} converges to the real solution $x^{\dag}$ of inverse problem (\ref{eq:1}).\\
(b) Consider that $ \phi^{T}$ cannot obtain the minimum value of (\ref{penalty}). Suppose Assumption \ref{assump3}, \ref{assump4} and \ref{assump1} hold, the measurement is clean, and $f_{\phi^{T}}(\cdot)$ is strongly convex. The sequence generated by Algorithm \ref{alg:2} converges to a $f_{\phi^{T}}$-minimizing solution of inverse problem (\ref{eq:1}).
\end{theorem}
\begin{remark}The above theorem proves the convergence of the proposed algorithm as the iterations converge to infinity. However, in practice, the number of iterations (number of layers of network unfolding) $K$ is limited. We will show the fast convergence of the Algorithm \ref{alg:2} later by numerical experiments, thus demonstrating that the algorithm output approximates the limit point well, even when $K$ is limited. \end{remark}

\begin{remark}In Theorem \ref{thm:2} (b), the strong convexity can be achieved by adding a term $\nu\|\cdot\|^2$ to $f_{\phi^{T}}(\cdot)$, where $\nu$ is a positive constant.\end{remark}

When the measurement contains noise, i.e., $\delta>0$, Proposition \ref{pro:1} shows that Algorithm \ref{alg:2} will perform $k_*$ iterations and then stops.
When the measurement is clean, Theorem \ref{thm:2} shows that Algorithm \ref{alg:2} will converge to a real solution or $f_{\phi^{T}}$-minimizing solution of inverse problem (\ref{eq:1}). Next, we will show that the $k_*$-step terminated solution is a stable approximation of the real solution or $f_{\phi^{T}}$-minimizing solution.
\begin{theorem}\label{thm:3}
Let the assumptions of Theorem \ref{thm:2} hold and let $k_*=k_*(\delta,y^\delta)$ be chosen according to the termination criterion (\ref{terminate}). Then, the output $x_{k_*}^\delta$ converges to the real solution or $ f_{\phi^{T}}$-minimizing solution of inverse problem (\ref{eq:1}) as $\delta\rightarrow0$.
\end{theorem}
\section{Implementation}\label{sect4}
The evaluation was performed on two multichannel MRI data with various trajectories. The details of the MRI data are as follows:
\subsection{Data Acquisition}
\subsubsection{Knee data}
First, we tested our proposed method on knee MRI data \footnote{\url{http://mridata.org/}}. The raw data were acquired from a 3T Siemens scanner. The number of coils was 15 and the 2D Cartesian turbo spin echo (TSE) protocol was used. The parameters for data acquisition are as follows: the repetition time (TR) was 2800ms, the echo time (TE) was 22ms, the matrix size was $768\times 770\times 1$ and the field of view (FOV) was $280 \times 280.7 \times 4.5 \text{mm}^3$. Particularly, the readout oversampling was removed by transforming the $k$-space to image, and cropping the center $384 \times 384$ region. Fully sampled multichannel knee images of nine volunteers were collected, of which data from seven subjects (including 227 slices) were used for training, while
the data from the remaining two subjects (including 64 slices) were used for testing.

\subsubsection{Human brain data}
To verify the generalization of the proposed method, we tested it on human brain MRI data \footnote{\url{https://drive.google.com/file/d/1qp-l9kJbRfQU1W5wCjOQZi7I3T6jwA37/view?usp=sharing}}, which was collected by \cite{8434321}. These MRI data were acquired using a
3D T2 fast spin echo with an extended echo train acquisition (CUBE) sequence with Cartesian readouts using a 12-channel head coil. The matrix dimensions were $256\times232\times208$ with 1 mm isotropic resolution.
The training data contain 360 slices $k$-space data from four subjects and the testing data contain 164 slices $k$-space data from two subjects. Each slice has a spatial dimension of $256\times232$.

\subsubsection{Sampling patterns}
\begin{figure}[thbp]
	\begin{center}
		\subfigure{\includegraphics[width=0.78\textwidth,height=0.25\textwidth]{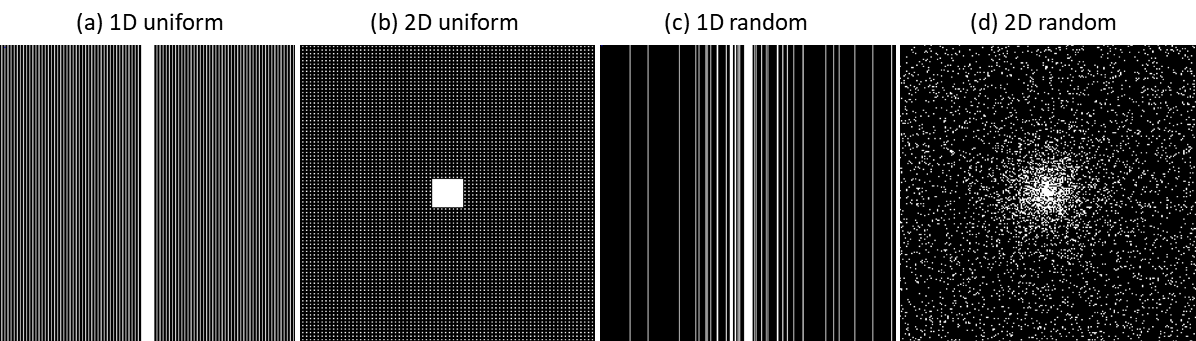}}
	\end{center}
	\caption{Various sampling patterns: (a) 1-D uniform undersampling at $R=4$, (b) 2-D uniform undersampling at $R=9$, (c) 1-D random undersampling at $R=8$, and (d) 2-D random undersampling at $R=12$.}
	\label{f2}
\end{figure}

Four types of undersampling patterns were considered, containing different acceleration factors and undersampling types. A visualization of these patterns is depicted in Figure \ref{f2}.
\subsection{Network Architecture and Training}
This paper considers the direct reconstruction of multichannel MR images to avoid additional calculations to estimate the coil sensitivities. The forward process of measurement acquisition can be modeled as the inverse problem:
$$PFx=y$$
where $x=[x_1,\ldots,x_C]$ denotes the $C$-channel MR image, $F$ denotes Fourier transform, and $P$ denotes the undersampling operator with above patterns. The unfolded PGD for above inverse problem reads:
 \begin{equation*}\begin{aligned}
\xi_{k+1}^n=&x_{k}^n-\eta F^{-1}P^*(PFx_{k}^n-y^n)\\
x_{k+1}^n=&\mathcal{S}_{\theta_{k+1}}(\xi_{k+1}^n).
\end{aligned}\end{equation*}
The above network is unfolded with ten layers (i.e., $K=10$), and the $\mathcal{S}_{\theta_{k+1}}$ module is represented by the following U-Net \cite{10.1007/978-3-319-24574-4_28}. Its encoder consists of $3 \times 3$ convolution, batch normalization \cite{pmlr-v37-ioffe15}, and ReLU and decoder part used the nearest neighbour upsampling and $3 \times 3$ convolution. The max pooling layers and skip connection through the channel concatenation are also included. We concatenate real and imaginary values along the channel dimension to handle complex values. Consequently, an input and an output have 2$C$ channels.
\begin{figure}[thbp]\label{unet}
	\begin{center}
		\subfigure{\includegraphics[width=0.55\textwidth,height=0.30\textwidth]{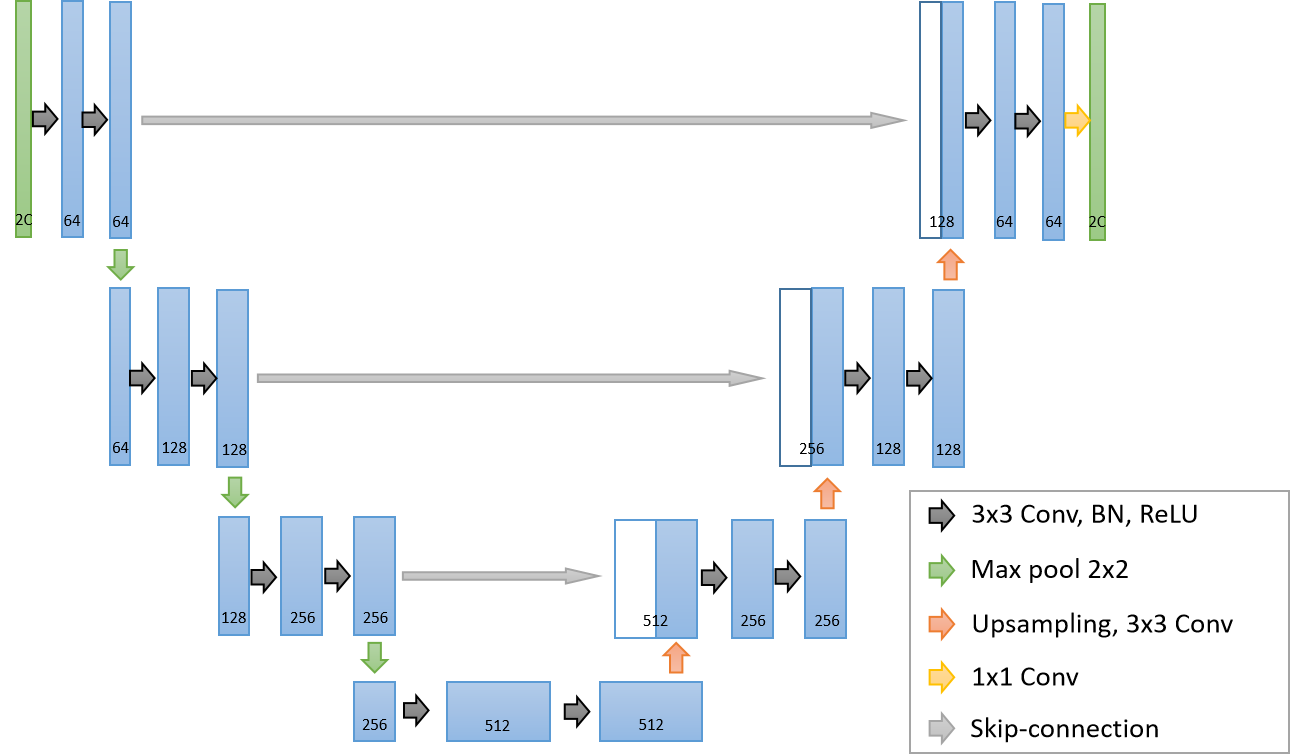}}
	\end{center}
	\caption{Network architecture for $\mathcal{S}_{\theta_{k+1}}$ module in unfolded PGD.}
\end{figure}

We designed the following ICNN to represent penalty function $f_{\phi}$. Its architecture combines PatchGAN discriminator \cite{Isola_2017_CVPR} and ICNN. The network has general and non-negative convolution layers, and batch normalization, ReLU and average-pooling.
\begin{figure}[thbp]
	\begin{center}
		\subfigure{\includegraphics[width=0.48\textwidth,height=0.28\textwidth]{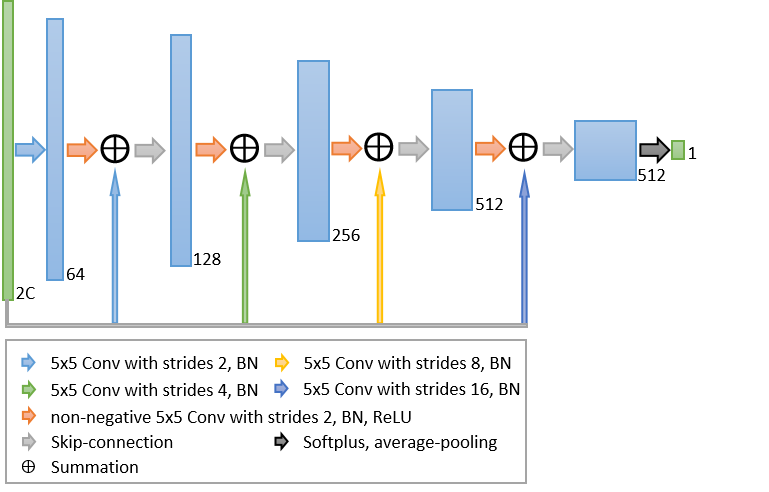}}
	\end{center}
	\caption{Network architecture for penalty function $f_{\phi}$.}
\end{figure}

In the training algorithm, we set $T_\Theta=2$ and $T_\phi=6$. For loss function $J_{1,t}$, we add a term $\mu_1\|\mathcal{S}_{\theta_k}(\xi_{k}^n)-x^n\|^2$ to increase the stability of the training. For loss function $J_{2,t+1}$, we add a gradient penalty term to constrain the 1-Lipschitz continuous of $f_{\phi}$. Namely,
\begin{equation*}\begin{aligned}J_{1,t}(\Theta)=&\sum_{n=1}^N\sum_{k=1}^K\left[\frac{1}{2}\|\mathcal{S}_{\theta_k}(\xi_{k}^n)-\xi_{k}^n\|^2+f_{\phi^t}(\mathcal{S}_{\theta_k}(\xi_{k}^n))+\mu_1\|\mathcal{S}_{\theta_k}(\xi_{k}^n)-x^n\|^2\right]\\
J_{2,t+1}(\phi)=&\mathbb{E}_{x^n\sim \mathbb{P}_x}[f_{\phi}(x^n)]-\frac{1}{K}\sum_{k=1}^K\left[\mathbb{E}_{x^n_0\sim \mathbb{P}_0}\left[f_\phi(\mathcal{S}_{\theta_k^{t+1}}(\xi^n_k))\right]+\mu_2\mathbb{E}_{\widehat{x}^k\sim \mathbb{P}_{\widehat{x}}^k}[(\|\nabla f_{\phi}(\widehat{x}^k)\|-1)^2]\right]
\end{aligned}\end{equation*}
where $\mathbb{P}_{\widehat{x}}^k$ sampling uniformly along straight lines between pairs of points sampled from the data distribution $\mathbb{P}_x$ and the distribution $\mathbb{P}_k$.

The ADAM \cite{kingma2014adam} optimizer with is chosen for Algorithm \ref{alg:1} with respect to the loss functions $J_{1,t}$ and $J_{2,t+1}$. The size of the mini batch is 1, and the number of epochs is 500. The learning rate $\gamma$ is set to $10^{-4 }$ and $\beta_1=0.9, \beta_2=0.999$.
The labels for the network were the fully sampled MR image. The input data for the network was the regridded downsampled k-space measurement from 1-D and 2-D random trajectories. Without specific instructions, we train the network separately for different trajectories.
The models were implemented on an Ubuntu 20.04 operating system equipped with an NVIDIA A6000 Tensor Core (GPU, 48 GB memory) in the open PyTorch 1.10 framework \cite{paszke2019pytorch} with CUDA 11.3 and CUDNN support.

\subsection{Performance Evaluation}
In this study, the quantitative evaluations were all calculated on the image domain. The image is derived using an inverse Fourier transform followed by an elementwise square-root of sum-of-the squares (SSoS) operation, i.e. $z[n]=(\sum_{i=1}^{C}|x_i[n]|^2)^{\frac{1}{2}}$, where $z[n]$ denotes the $n$-th element of image $z$, and $x_i[n]$ denotes the $n$-th element of the $i$th coil image $x_i$. For quantitative evaluation, the peak signal-to-noise ratio (PSNR), normalized mean square error (NMSE) value and structural similarity (SSIM) index \cite{1284395} were adopted.
\section{Experimentation Results }\label{sect5}
\subsection{Comparative Studies}
In this section, we test our proposed new training, and testing procedures for PGD unfolded network, dubbed PGD-Net+ for short, on the knee and brain datasets. A series of extensive comparative experiments were studied to demonstrate our methods' effectiveness. In particular, we compared the traditional parallel imaging (PI) method (ESPIRiT \cite{uecker2014espirit}), and SOTA unfolded network, i.e., MoDL \cite{8434321}, where the network module of MoDL was replaced with U-Net architecture, as shown in Figure \ref{unet}, for fairness of comparison. In addition, because the proposed method is formally like the GAN model, i.e., the learned penalty function is similar to a discriminator, and the unfolding network is similar to a generator, this paper also compares it with the CycleGAN model. The CycleGAN uses a U-Net for the generator, and a pyramid network for the discriminator \cite{oh2020unpaired}. Finally, to verify the effectiveness of the introduced penalty function, we used a PGD unfolded network trained with $L_2$-norm (\ref{loss2}) as the loss function, dubbed PGD-Net, as an ablation experiment.
\subsection{Experiments on Clean Measurement}
\begin{figure}[!t]
\centering
\includegraphics[width=0.98\textwidth,height=0.5\textwidth]{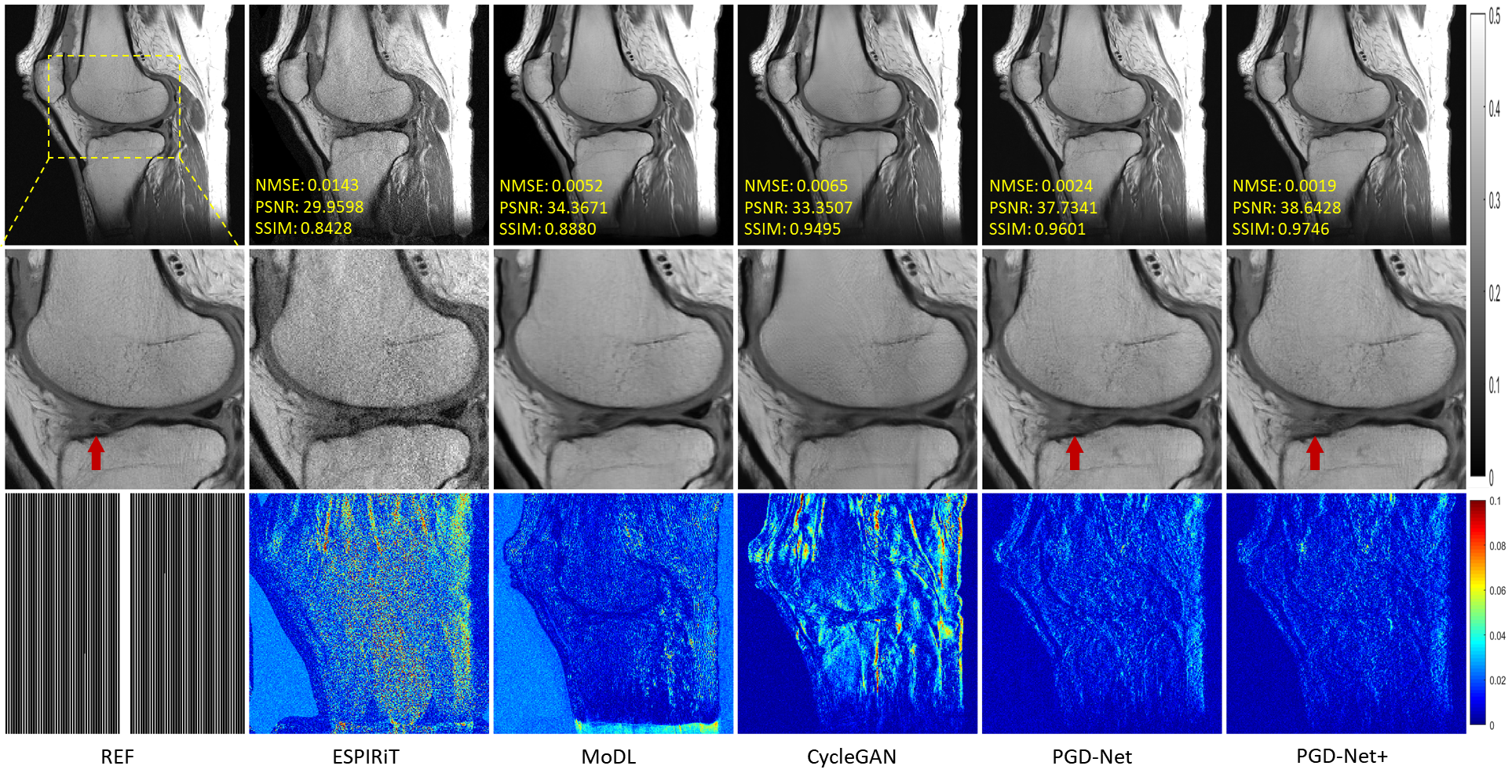}
\caption{Reconstruction results under 1-D uniform undersampling at $R=4$. The values in the corner are the NMSE/PSNR/SSIM values of each slice. The second and third rows illustrate the enlarged and error views, respectively. The grayscale of the reconstructed images and the color bar of the error images are at the right of the figure.}
\label{f4}
\end{figure}

This section tests our proposed PGD-Net+ and comparative networks or algorithms on clean measurement.
Figure \ref{f4} shows the reconstruction results of the knee data using various methods under a 1-D uniform trajectory with an acceleration factor of 4. As shown in Figure \ref{f4}, the aliasing pattern remains in the reconstructed images for the traditional ESPIRiT method. As seen from the enlargement, the CycleGAN and MoDL reconstructions are too smooth, and some image detail textures are lost. For the ablation experiment, although the quality of the PGD-Net reconstruction is close to that of our proposed method, the detail of our method is more clearly recovered, as can be seen in the knee cartilage indicated by the arrow.

\begin{figure}[!t]
\centering
\includegraphics[width=0.98\textwidth,height=0.5\textwidth]{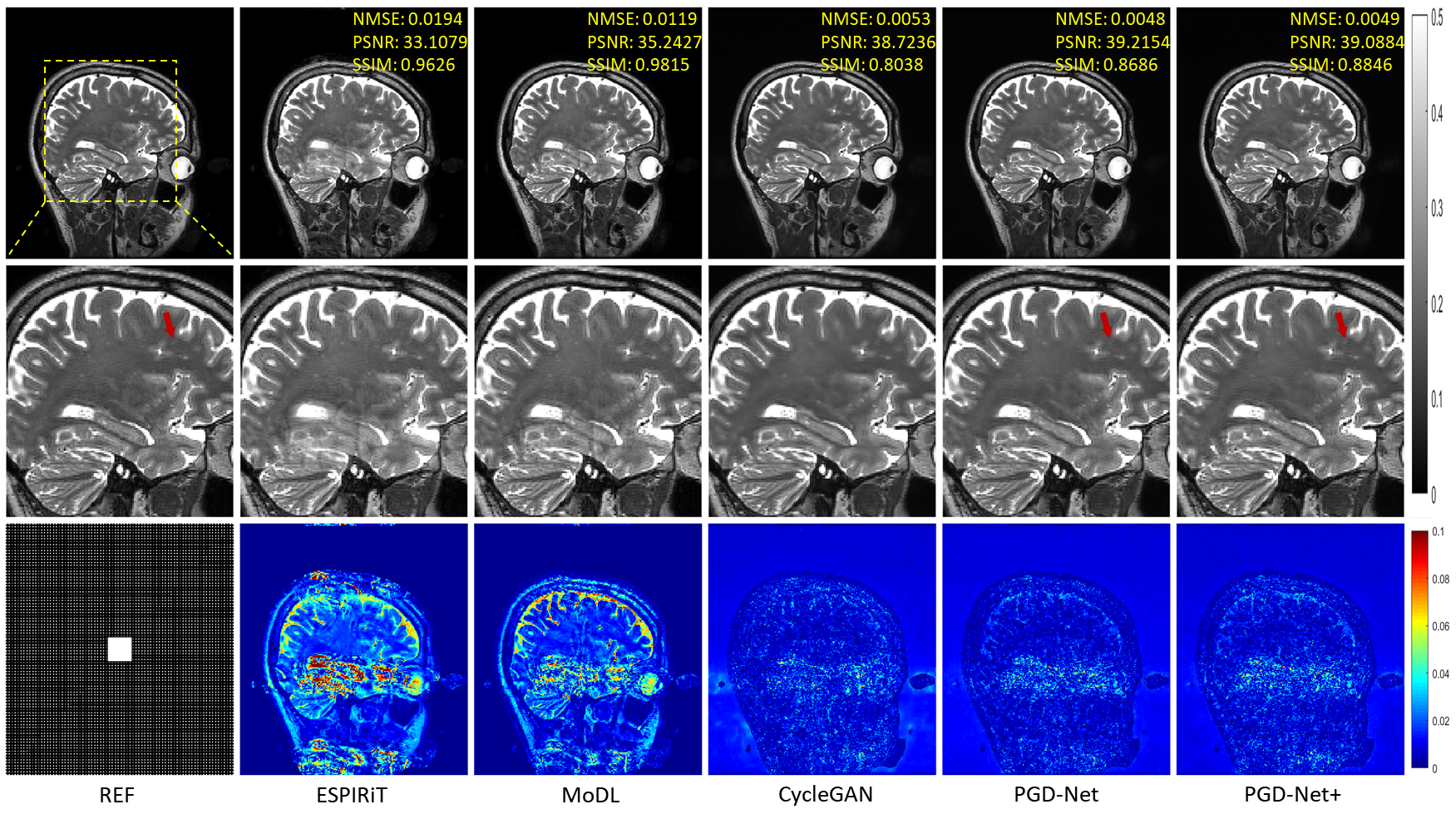}
\caption{Reconstruction results under 2-D uniform undersampling at $R=9$. The values in the corner are the NMSE/PSNR/SSIM values of each slice. The second and third rows illustrate the enlarged and error views, respectively. The grayscale of the reconstructed images and the color bar of the error images are at the right of the figure.}
\label{f5}
\end{figure}

Figure \ref{f5} shows the reconstruction results of the brain data using various methods under a 2-D uniform trajectory with an acceleration factor of 9. It is easy to see that the aliasing pattern remains in the ESPIRiT and MoDL reconstructions, while the CycleGAN reconstruction is slightly smooth. In terms of reconstruction error, the remaining two methods reconstruct the image well. However, the enlarged image makes it easy to see that the proposed method reconstructs the details at the point indicated by the arrow more accurately.

The competitive quantitative results of the above methods are shown in Table \ref{tab:1}. Our method consistently outperforms the other comparative methods for knee and brain data. Therefore, as characterized by visual and quantitative evaluations, the above experiments confirm the competitiveness of our method under the case with no additional noise on measurement.
\begin{table}\label{tab:1}
	\begin{center}
		\caption{Quantitative comparison for various methods on the knee and brain data.}\label{tab:1}
		\setlength{\tabcolsep}{3mm}{
			\begin{tabular}{l|l|ccc}
				\hline
				\multicolumn{ 2}{c}{ Datasets} & \multicolumn{ 3}{|c}{Quantitative Evaluation}  \\
				\multicolumn{ 2}{c|}{ \& Methods   } &NMSE &PSNR(dB)&SSIM  \\
				\hline
				\multirow{5}{*}{Knee}
				& ESPIRiT  &0.0119$\pm$0.0053&31.53$\pm$1.94&0.79$\pm$0.09. \\
				\cline{2-5}
				& MoDL  &0.0071$\pm$0.0039&33.71$\pm$1.54&0.80$\pm$0.11\\
				\cline{2-5}
				& CycleGAN &0.0086$\pm$0.0041&33.00$\pm$2.38&0.94$\pm$0.02\\
				\cline{2-5}
				&  PGD-Net &0.0021$\pm$0.0010&39.07$\pm$1.96&0.94$\pm$0.03\\
				\cline{2-5}
				&  PGD-Net+ &\textcolor{red}{0.0020$\pm$0.0008}&\textcolor{red}{39.27$\pm$1.61}&0.93$\pm$0.04\\
				\hline
				\multirow{5}{*}{Brain}
				& ESPIRiT  &0.0121$\pm$0.0050&35.61$\pm$3.15&0.96$\pm$0.02 \\
				\cline{2-5}
				& MoDL  &0.0048$\pm$0.0021&39.32$\pm$1.89&0.97$\pm$0.06\\
				\cline{2-5}
				& CycleGAN &0.0062$\pm$0.0036&38.24$\pm$1.26&0.82$\pm$0.04\\
				\cline{2-5}
				&  PGD-Net &0.0039$\pm$0.0013&40.19$\pm$1.27&0.88$\pm$0.02\\
				\cline{2-5}
				& PGD-Net+ &\textcolor{red}{0.0038$\pm$0.0014}&\textcolor{red}{40.30$\pm$1.49}&0.89$\pm$0.02\\
				\hline
		\end{tabular}}
	\end{center}
\end{table}

\begin{figure}[!t]
\centering
\includegraphics[width=0.98\textwidth,height=0.5\textwidth]{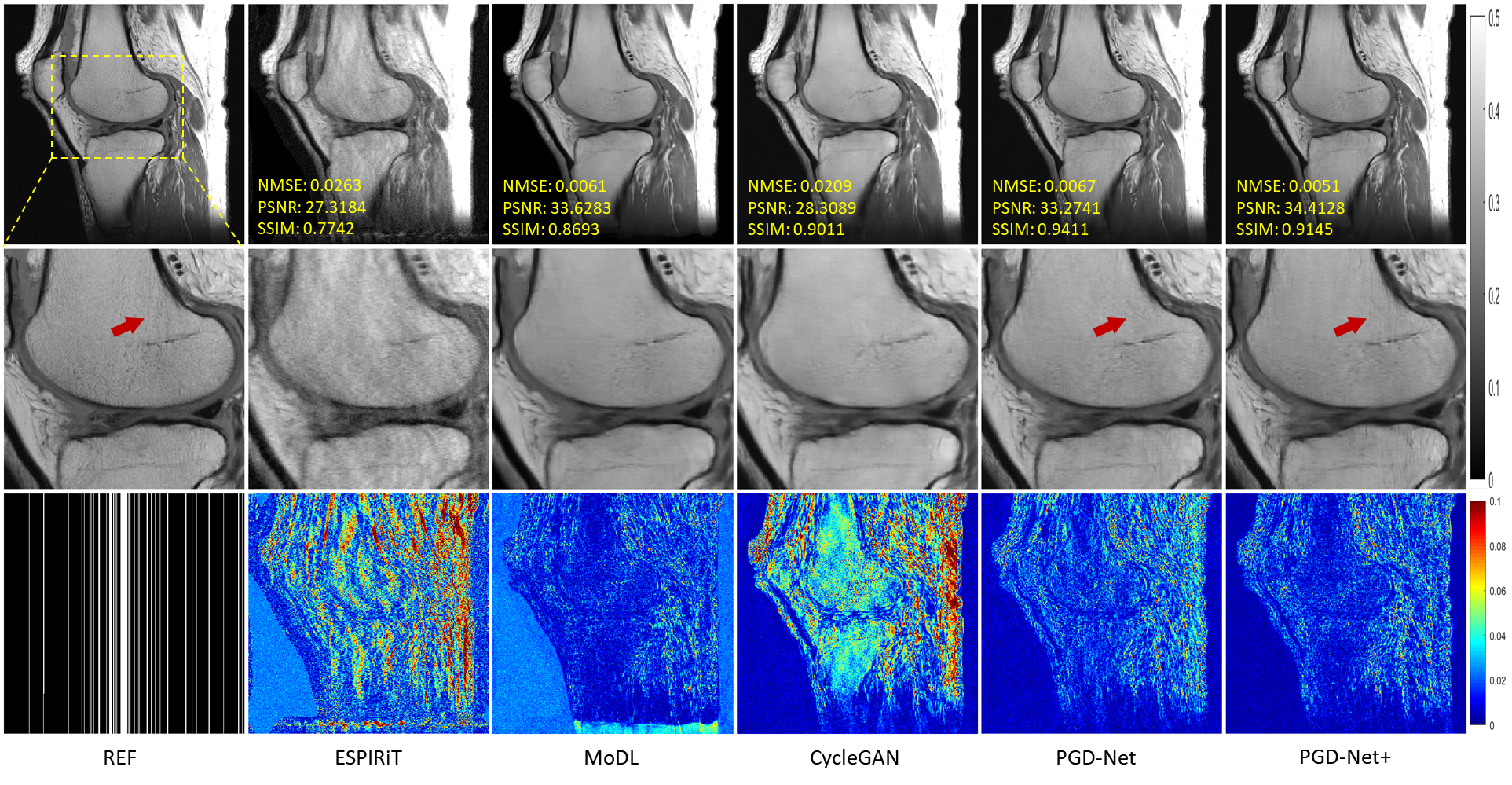}
\caption{Reconstruction results under 1-D random undersampling at $R=8$. The values in the corner are the NMSE/PSNR/SSIM values of each slice. The second and third rows illustrate the enlarged and error views, respectively. The grayscale of the reconstructed images and the color bar of the error images are at the right of the figure.}
\label{f6}
\end{figure}

\begin{figure}[!t]
\centering
\includegraphics[width=0.98\textwidth,height=0.5\textwidth]{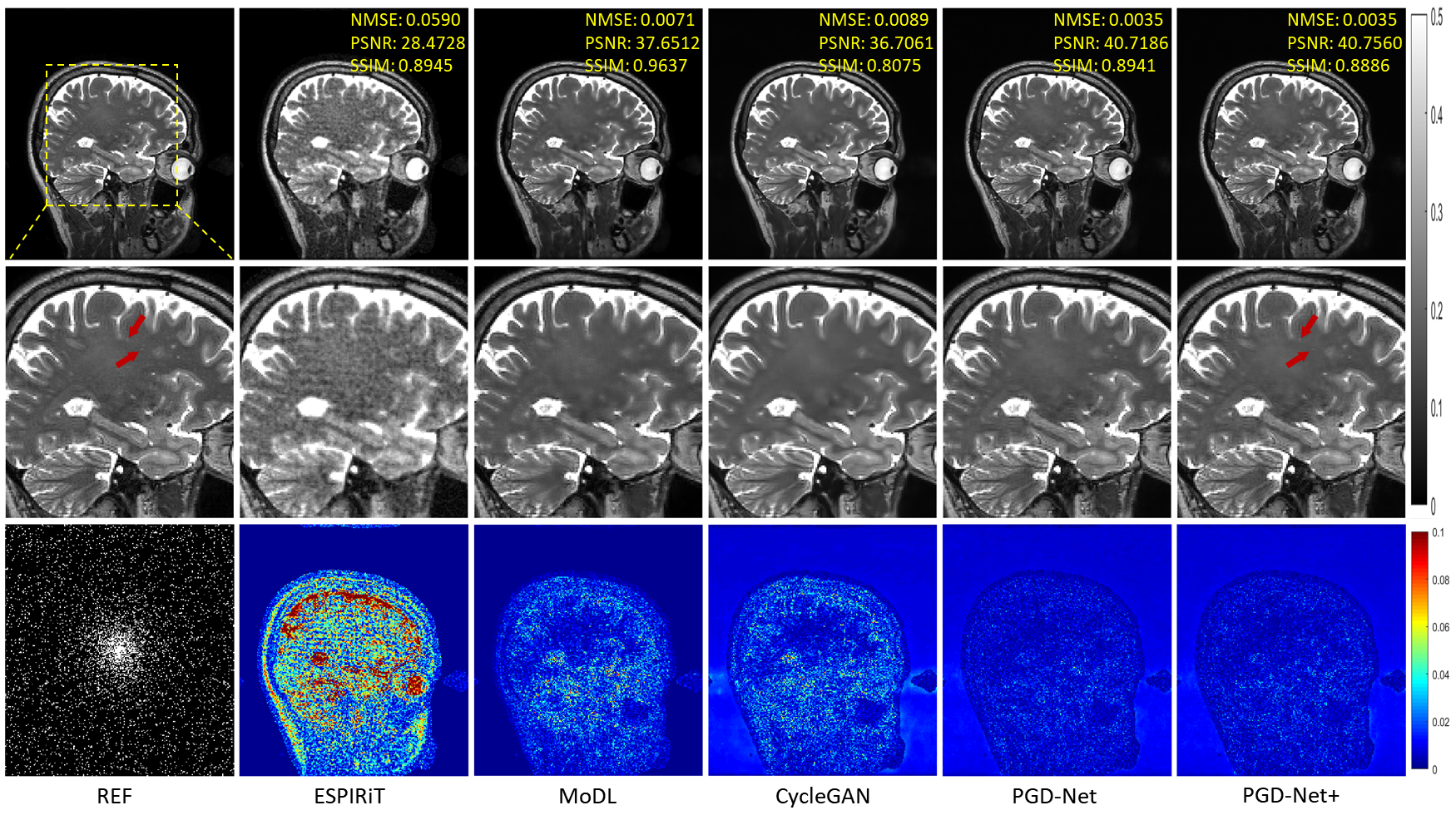}
\caption{Reconstruction results under 2-D random undersampling at $R=12$. The values in the corner are the NMSE/PSNR/SSIM values of each slice. The second and third rows illustrate the enlarged and error views, respectively. The grayscale of the reconstructed images and the color bar of the error images are at the right of the figure.}
\label{f7}
\end{figure}

To further validate the effectiveness of the proposed method, we tested the reconstruction results at higher acceleration rates. Figure \ref{f6} shows the reconstruction results of the knee data using various methods under a 1-D random trajectory with an acceleration factor of 8. Similar to the results of the previous experiment, the ESPIRiT reconstruction still leaves the aliasing pattern, and the MoDL and CycleGAN reconstructions are too smooth. PGD-Net and the proposed method reconstruct well, but the proposed method performs better in the texture reconstruction indicated by the arrow. Higher acceleration rates were also tested on the brain data for the reconstruction. Figure \ref{f7} shows the reconstruction results of the brain data using various methods under a 2-D random trajectory with an acceleration factor of 12. It is not difficult to find that of all the methods, ours is the only one that reconstructs the detail indicated by the arrow.

Thus far, we have confirmed the competitiveness of the proposed method in terms of reconstruction quality by experimenting with different acceleration rates and sampling patterns.
\subsection{Convergence Speed}
Above, we claimed that "the proposed method achieves faster convergence because the loss function measures the distance from each iteration to the real data manifold, rather than just the $L_2$-norm between the last iteration and the real data". Then, we will design experiments to verify it.

\begin{figure}[!t]
\centering
\includegraphics[width=0.98\textwidth,height=0.35\textwidth]{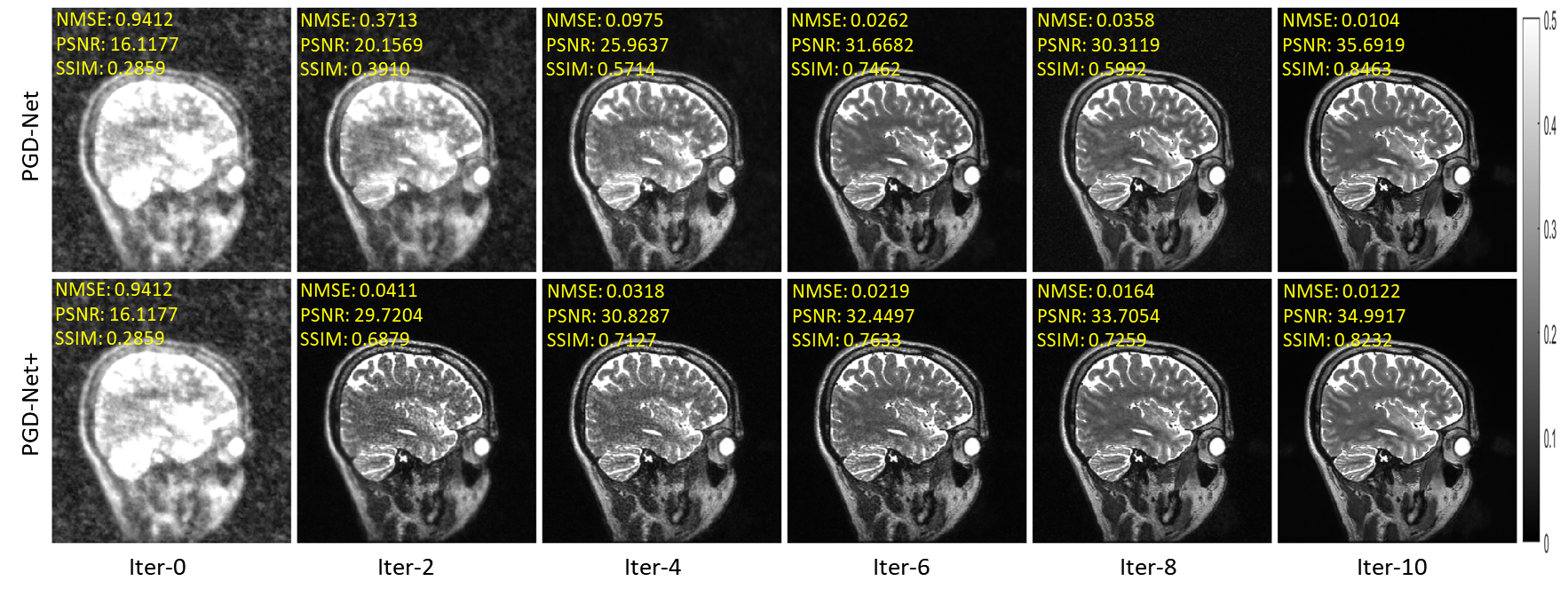}
\caption{Different iterations of reconstructions under 2-D random undersampling at $R=12$. The values in the corner are each slice's NMSE/PSNR/SSIM values. The grayscale of the reconstructed images is at the right of the figure.}
\label{f8}
\end{figure}

\begin{figure}[!t]
\centering
\includegraphics[width=0.98\textwidth,height=0.35\textwidth]{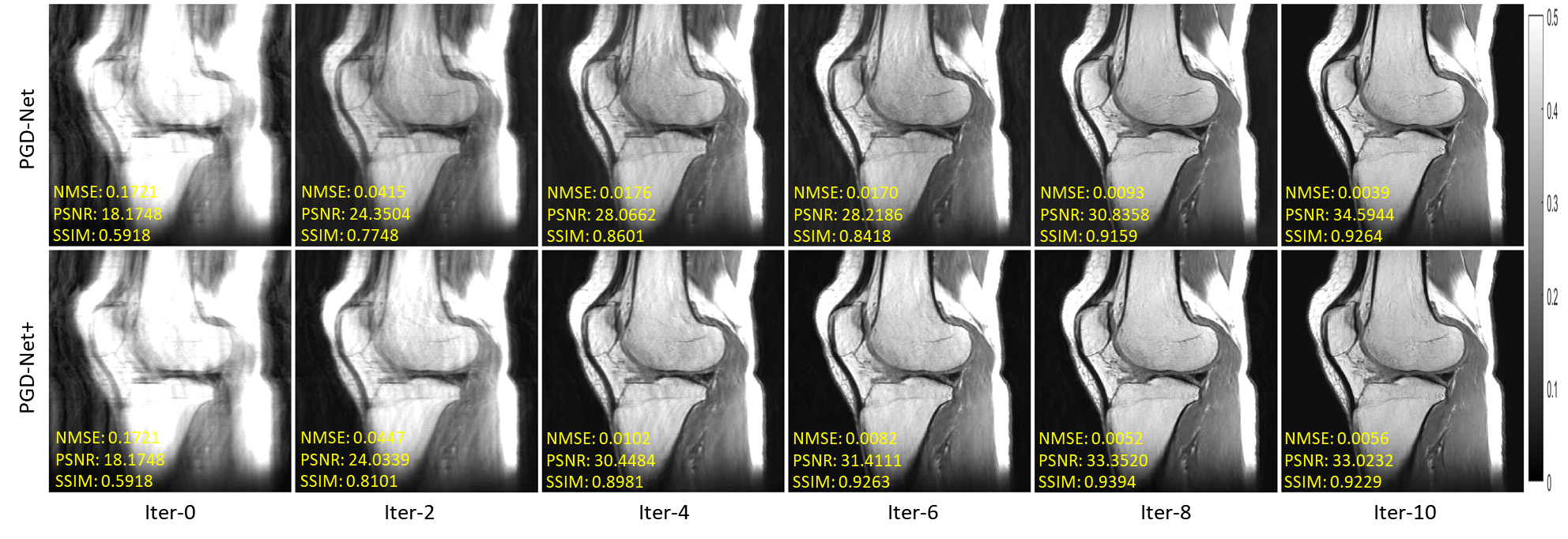}
\caption{Different iterations of reconstructions under 1-D random undersampling at $R=8$. The values in the corner are each slice's NMSE/PSNR/SSIM values. The grayscale of the reconstructed images is at the right of the figure..}
\label{f9}
\end{figure}

\begin{figure}[!t]
\centering
\includegraphics[width=0.78\textwidth,height=0.35\textwidth]{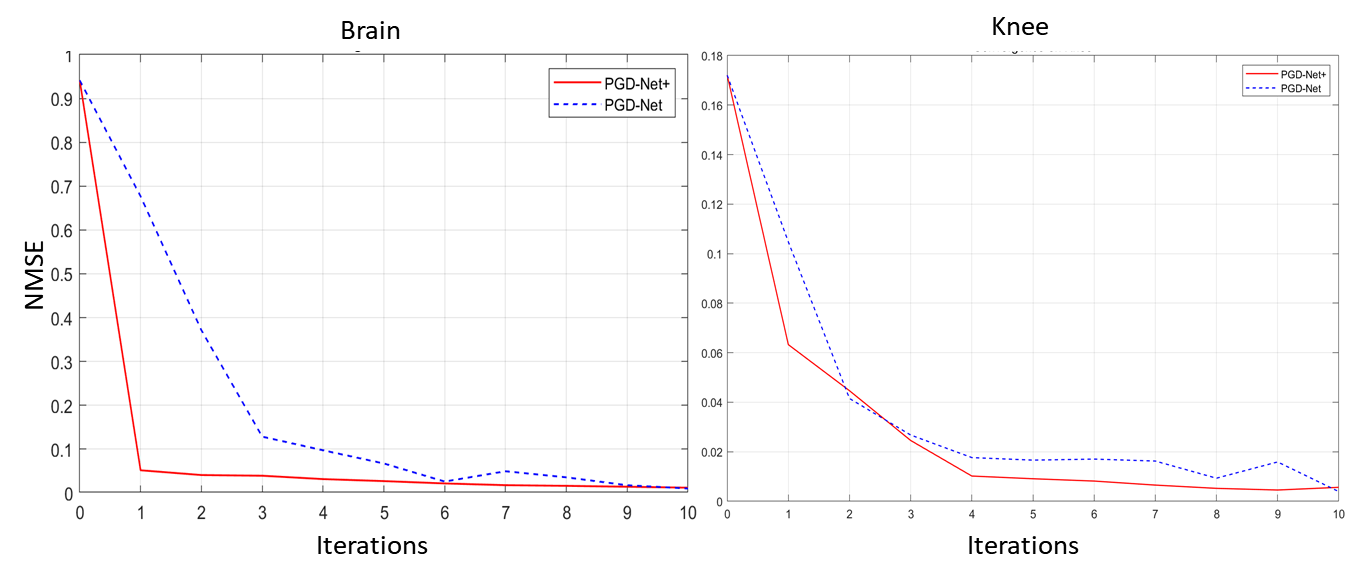}
\caption{Illustration of NMSE decreases with iteration executed for both methods on brain and knee data when the measurements are clean.}
\label{f10}
\end{figure}

Figure \ref{f8} shows the results of PGD-Net versus the proposed method for different iterations on 12-fold accelerated brain data. It is easy to see that the second iteration of our method gives a clear reconstruction, while the second iteration of PGD-Net still gives a very blurred reconstruction. It is also evident from the quantitative metrics that the proposed method pulls away from PGD-Net from the second iteration onwards. Figure \ref{f9} shows the results of the two methods for different iterations on 8-fold accelerated knee data. It is also easy to see from Figure \ref{f9} that the aliasing pattern is well suppressed in the fourth iteration of the proposed method, while the aliasing pattern remains evident in the fourth iteration of PGD-Net. In addition, Figure \ref{f10} shows how the NMSE decreases with iteration executed for both methods on brain and knee data. It is clear that the proposed method converges faster than PGD-Net.

\subsection{Experiments on Noisy Measurement}
If the measurement contains noise, the iterative algorithm usually has semi-convergence when solving the inverse problem, i.e., it will start the iteration by travelling towards the true solution but will deviate from it as the iteration progresses. Therefore, a suitable termination criterion is important.

\begin{figure}[!t]
\centering
\includegraphics[width=0.98\textwidth,height=0.35\textwidth]{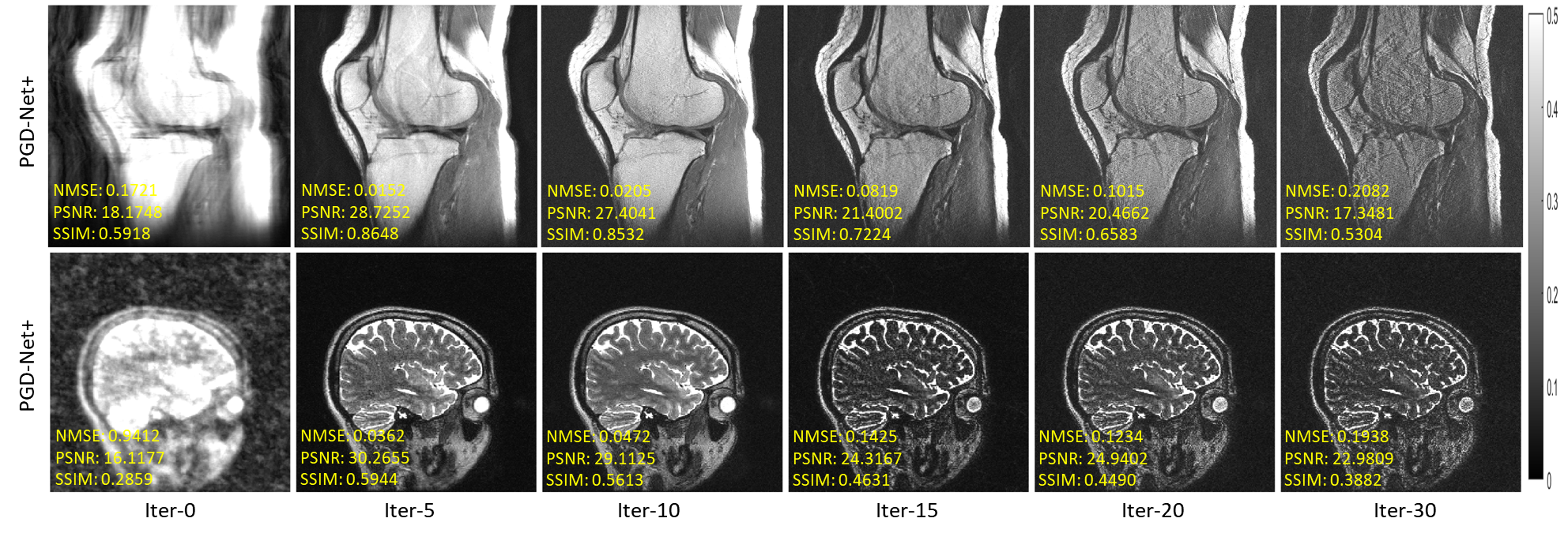}
\caption{Different iterations of reconstructions under 8-fold 1-D and 12-fold 2-D random undersampling, respectively, when the measurement is noisy. The values in the corner are each slice's NMSE/PSNR/SSIM values. The grayscale of the reconstructed images is at the right of the figure.}
\label{f11}
\end{figure}

\begin{figure}[!t]
\centering
\includegraphics[width=0.78\textwidth,height=0.35\textwidth]{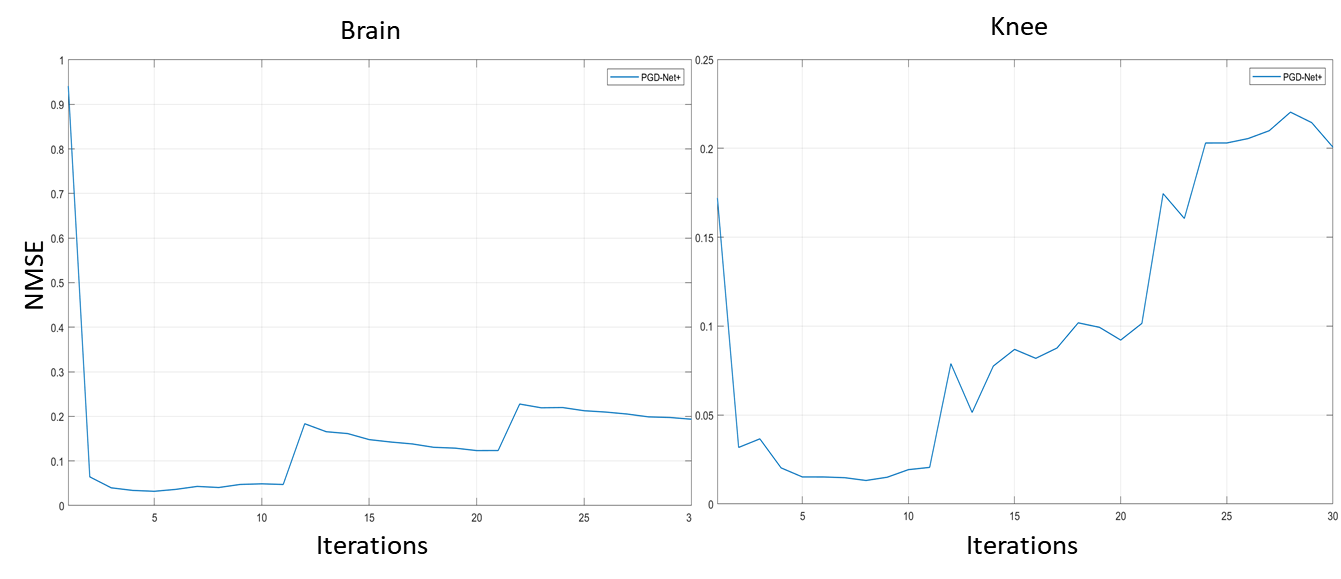}
\caption{Illustration of the walk of NMSE with iteration executed for both methods on brain and knee data when the measurements are noisy.}
\label{f12}
\end{figure}

Figure \ref{f11} shows the performance of our PGD ($K=10$ at training) after performing more iterations when the measurement contains $2.5\%$ noise (i.e., $\|n\|/\|y\|=0.025$). In terms of visual perception, the reconstruction of PGD-Net+ is optimal for ten iterations on the knee data and five iterations on the brain data, after which the reconstruction quality gradually decreases.
Figure \ref{f12} shows the NMSE curves with iteration executed for both methods on brain and knee data. The NMSE curve is consistent with the visual perception of the reconstructed image.

In Algorithm \ref{alg:2}, if $\tau$ is chosen to be 2 on the knee data and $\tau$ is chosen to be 2.5 on the brain data, Algorithm \ref{alg:2} terminates at step 9 on the knee data and terminates at step 4 on the brain data. The results of this experiment provide a good validation of Theorem \ref{thm:3} and show that the proposed method can produce an acceptable reconstruction result even when the measurement data contains noise.
\section{Conclusion}\label{sect6}
In this paper, we proposed a new training procedure for the unfolding method and proved it is an iterative regularization method. The working mechanism of the proposed method is interpretable. Namely, the penalty function represented by ICNN is used to measure the distance to the real data manifold, and the network module in the unfolding method is used to approximate the proximal mapping of this penalty function. Suppose the real data manifold intersects the inverse problem solutions with only the unique real solution. We proved that the unfolded method would converge to it stably. Experimentally, we demonstrated with an example of MRI reconstruction that the proposed method outperforms conventional unfolding methods and traditional regularization methods in terms of reconstruction quality, stability and convergence speed.

However, there are at least two points of future work for us as follows. First, in theory, the proposed method wants a unique real solution to the inverse problem relying on strong Assumption \ref{assump2}, which is difficult to verify in practice. In compressive sensing theory, many conditions guarantee the uniqueness of the reconstruction that has been proposed, such as RIP, NSP, etc. Therefore, it will be one of our future works to construct new mechanisms for the DL method based on these conditions to guarantee the uniqueness of solving the inverse problem. Furthermore, experimentally, the effectiveness of the proposed method is verified in this paper by MRI only. Extending this method to other imaging modalities, such as CT and PET, e.t.c., is also our future work.
\section*{Appendix}
\subsection*{The proof of Theorem \ref{thm:1}}
\begin{proof}
This proof follows the technical route in \cite{NEURIPS2018_d903e960}. First, we prove that the distance function $d_{\mathcal{M}}$ is convex and 1-Lipschitz continuous. By the definition of $d_{\mathcal{M}}$, for any $\lambda\in(0,1)$, $x,y\in X$, we have
\begin{equation*}\begin{aligned}
d_{\mathcal{M}}(\lambda x+(1-\lambda)y)&=\min_{z\in\mathcal{M}}\|\lambda x+(1-\lambda)y-z\|\\
&\leq\|\lambda x+(1-\lambda)y-\lambda x^*+(1-\lambda)y^*\|\\
&\leq\lambda\|x-x^*\|+(1-\lambda)\|y-y^*\|\\
&=\lambda d_{\mathcal{M}}(x)+(1-\lambda)d_{\mathcal{M}}(y)
\end{aligned}\end{equation*}
where $x^*\in\arg\min_{z\in \mathcal{M}}\|z-x\|$ and $y^*\in\arg\min_{z\in \mathcal{M}}\|z-y\|$. For any $x,y\in X$, we have
\begin{equation*}\begin{aligned}
d_{\mathcal{M}}(x)-d_{\mathcal{M}}(y)=&\min_{z\in\mathcal{M}}\|x-z\|-\min_{z\in\mathcal{M}}\|y-z\|\\
\leq&\|x-y^*\|-\|y-y^*\|\leq\|x-y\|
\end{aligned}\end{equation*}
where $y^*\in\arg\min_{z\in \mathcal{M}}\|z-y\|$. Thus, we have proved the convexity and 1-Lipschitz continuity of $d_{\mathcal{M}}$.

Next, we prove that the minimum of functional (\ref{adv}) is $d_{\mathcal{M}}$. For any convex and 1-Lipschitz continuous non-negative function $f$, we have
\begin{equation*}\begin{aligned}
\frac{1}{K}\sum_{k=1}^K\mathbb{E}_{Z^k\sim \mathbb{P}_k}\left[f(Z^k)\right]-\mathbb{E}_{X\sim \mathbb{P}_x}\left[f(X)\right]=&\frac{1}{K}\sum_{k=1}^K\mathbb{E}_{Z^k\sim \mathbb{P}_k}\left[f(Z^k)-f(P_{\mathcal{M}}(Z^k))\right]\\
\leq&\frac{1}{K}\sum_{k=1}^K\mathbb{E}_{Z^k\sim \mathbb{P}_k}\left[\|Z^k-P_{\mathcal{M}}(Z^k)\|\right]\\
=&\frac{1}{K}\sum_{k=1}^K\mathbb{E}_{Z^k\sim \mathbb{P}_k}\left[d_{\mathcal{M}}(Z^k)\right]\\
=&\frac{1}{K}\sum_{k=1}^K\mathbb{E}_{Z^k\sim \mathbb{P}_k}\left[d_{\mathcal{M}}(Z^k)\right]-\mathbb{E}_{X\sim \mathbb{P}_x}\left[d_{\mathcal{M}}(X)\right]
\end{aligned}\end{equation*}
where the first equality is due to Assumption (\ref{assump4}), the first inequality is due to the 1-Lipschitz continuity of $f$, the second equality is due to the definition of $d_{\mathcal{M}}$ and the last equality is due to the fact that $\mathbb{E}_{X\sim \mathbb{P}_x}\left[d_{\mathcal{M}}(X)\right]=0$. The above inequality shows that $d_{\mathcal{M}}$ is the convex and 1-Lipschitz continuous non-negative minimum of functional (\ref{adv}).
\end{proof}
\subsection*{The proof of Proposition \ref{pro:1}}
\begin{proof}
Let $x^*$ be the $ f_{\phi^{T}}$-minimizing solution of inverse problem (\ref{eq:1}), we have
\begin{equation}\label{eq:app:1}\begin{aligned}
&\frac{1}{2}\|x_{k+1}^\delta-x^*\|^2-\frac{1}{2}\|x_k^\delta-x^*\|^2\\
=&\frac{1}{2}\|x_{k+1}^\delta-x^*\|^2-\frac{1}{2}\|x_{k}^\delta-x_{k+1}^\delta\|^2-\frac{1}{2}\|x_{k+1}^\delta-x^*\|^2-\langle x_k^\delta-x_{k+1}^\delta,x_{k+1}^\delta-x^*\rangle\\
=&-\frac{1}{2}\|x_{k}^\delta-x_{k+1}^\delta\|^2-\langle\eta A^*(Ax_k-y^\delta)+\xi_{k+1},x_{k+1}^\delta-x^* \rangle\\
\leq&-\frac{1}{2}\|x_{k}^\delta-x_{k+1}^\delta\|^2+f_{\phi^{T}}(x)-f_{\phi^{T}}(x_{k+1}^\delta)+\epsilon_k-\eta\langle A^*(Ax_k^\delta-y^\delta),x_{k+1}^\delta-x^* \rangle\\
\leq&-\frac{1}{2}\|Ax_k^\delta-Ax_{k+1}^\delta\|^2+f_{\phi^{T}}(x)-f_{\phi^{T}}(x_{k+1}^\delta)+\epsilon_k-\eta\|Ax_k^\delta-y^\delta\|^2\\
&-\eta\langle Ax_k^\delta-y^\delta,Ax_{k+1}^\delta-Ax_k^\delta+y^\delta-y \rangle\\
\leq&f_{\phi^{T}}(x^*)-f_{\phi^{T}}(x_{k+1}^\delta)-\left(\eta-\frac{\eta^2}{2}\right)\|Ax_k^\delta-y^\delta\|^2+\eta\delta\|Ax_k^\delta-y^\delta\|+\epsilon_k\\
\end{aligned}\end{equation}
where the first inequality is due to Assumption \ref{assump1}, the second inequality is due to $\|A\|\leq1$ and the last inequality is due to the fact that $-a^2-2\langle a,b \rangle\leq b^2$.

Assume that $k_*$ is the maximum number of iteration steps before Algorithm \ref{alg:2} triggers the termination criterion (\ref{terminate}), we have
\begin{equation*}
(k_*-1)\theta\tau^2\delta^2\leq\sum_{i=1}^{k_*-1}\left[\left(\eta-\frac{3\eta}{2\tau}-\frac{\eta^2}{2}\right)\|Ax_i^\delta-y^\delta\|^2+ \left(1-\frac{\eta}{2\tau}\right)(f_{\phi^{T}}(x_{i}^\delta)-f_{\phi^{T}}(x^*))\right] <\sum_{i=1}^{k_*}\epsilon_k+\frac{1}{2}\|x_0^\delta-x^*\|^2
\end{equation*}
where $\theta=\min\{1-2\tau/\eta,\eta-{\eta}/{\tau}-{\eta^2}/{2}\}$ the second inequality is due to $\|Ax_{i}^\delta-y^\delta\|+f_{\phi^{T}}(x_{i}^\delta)-f_{\phi^{T}}(x^*))>\tau\delta$ and the relation $2ab<a^2+b^2$, for any $i\leq k_*$.
\end{proof}
\subsection*{The proof of Theorem \ref{thm:2}}
\begin{proof}
When the measurement contains no noise, i.e., $\delta = 0$, from inequality (\ref{eq:app:1}), we have
\begin{equation*}
f_{\phi^{T}}(x_{k+1})-f_{\phi^{T}}(x^*)+\left(\eta-\frac{\eta^2}{2}\right)\|Ax_k-y\|^2\leq\frac{1}{2}\|x_k-x^*\|^2-\frac{1}{2}\|x_{k+1}-x^\dag\|^2+\epsilon_k\\
\end{equation*}
Summing both sides of the above inequality from 0 to $+\infty$, we get:
\begin{equation}\label{eq:app:3}\sum_{k=1}^{+\infty}\left[f_{\phi^{T}}(x_{k+1})-f_{\phi^{T}}(x^*)+\left(\eta-\frac{\eta^2}{2}\right)\|Ax_k-y\|^2\right]\leq\frac{1}{2}\|x_{0}-x^*\|^2+\sum_{k=1}^{+\infty}\epsilon_k\end{equation}
Then, we have
\begin{equation}\label{eq:app:2}\lim_{k\rightarrow+\infty}f_{\phi^{T}}(x_{k})=f_{\phi^{T}}(x^*)~\text{and}~\lim_{k\rightarrow+\infty}\|Ax_k-y\|=0\end{equation}

If $ \phi^{T}$  reaches the minimum of (\ref{penalty}) in Algorithm \ref{alg:1}, Theorem \ref{thm:1} shows that $f_{\phi^{T}}(x_k)=d_{\mathcal{M}}(x_k)$. Define $x'_k=P_{\mathcal{M}}(x_k)$ and $x''_k=\arg\min_{x\in\{x\in X| Ax=y \}}\|x_k-x\|$, we have $$\|x'_k-x''_k\|\leq \|x'_k-x_k\|+\|x''_k-x_k\|=d_{\mathcal{M}}(x_k)+\|x''_k-x_k\|$$
From (\ref{eq:app:2}) we know that $\lim_{k\rightarrow+\infty}\|x'_k-x''_k\|=0$. If further assumption \ref{assump2} holds, we have $\lim_{k\rightarrow+\infty}x'_k=\lim_{k\rightarrow+\infty}x''_k=x^*=x^\dag$. Then, we have
$$\|x_k-x^\dag\|\leq d_{\mathcal{M}}(x_k)+\|x'_k-x^{\dag}\|$$
which implies that $\|x_k-x^\dag\|\rightarrow0$ as $k\rightarrow+\infty$.

If $ \phi^{T}$ cannot reach the minimum  of (\ref{penalty}) in Algorithm \ref{alg:1} but $f_{\phi^{T}}$ is $\sigma$-strongly convex, inequality (\ref{eq:app:3}) implies that
$$\sum_{k=1}^{+\infty}\sigma\|x_k-x^*\|^2\leq\frac{1}{2}\|x_{0}-x^*\|^2+\sum_{k=1}^{+\infty}\epsilon_k. $$
Then, we show that $x_k$ converges to a $f_{\phi^{T}}$-minimizing solution of inverse problem (\ref{eq:1}) as $k\rightarrow+\infty$.
\end{proof}
\subsection*{The proof of Theorem \ref{thm:3}}
\begin{proof}
Let $\delta_n$ be a sequence converging to 0 as $n\rightarrow+\infty$ and $k_n=k_*(\delta_n,y^{\delta_n})$ be the maximum number of iteration steps before Algorithm \ref{alg:2} triggers the termination criterion (\ref{terminate}).

For any finite accumulation point $k$ of $\{k_n\}$, without loss of generality, we can assume that $k_n=k$ for all $n$. As $k$ fixed, $x^{\delta_n}_k$ depends continuously on $y^{\delta_n}$, we have
$$x^{\delta_n}_k\rightarrow x_k~\text{and}~\|Ax_{k}^{\delta_n}-y^{\delta_n}\|^2+f_{\phi^T}(x_{k}^{\delta_n})-f_{\phi^T}^*\rightarrow 0$$
as $n\rightarrow+\infty$, which means that $x_k$ is a $ f_{\phi^{T}}$-minimizing solution of inverse problem (\ref{eq:1}).

Consider the case $k_n\rightarrow+\infty$ as $n\rightarrow+\infty$. For $k<k_n$, we have
$$\|x_{k_n}^\delta-x^*\|\leq\|x_{k}^\delta-x^*\|+\sum_{i=k+1}^{k_n}\epsilon_i\leq\|x_{k}^\delta-x_k\|+\|x_k-x^*\|+\sum_{i=k+1}^{+\infty}\epsilon_i$$
where the first inequality is due to the inequality (\ref{eq:app:1}). Given $\varepsilon>0$, we can fix some $k_1(\varepsilon)$ such that $\sum_{i=k+1}^{+\infty}\epsilon_i<\varepsilon/3$ for any $k\geq k_1(\varepsilon)$, and we can also fix some $k_2(\varepsilon)$ such that $\|x_k-x^*\|<\varepsilon/3$ for any $k\geq k_2(\varepsilon)$. Because $x^{\delta_n}_k$ depends continuously on $y^{\delta_n}$, we can fix some $n(\varepsilon,k)$ such that $\|x_{k}^\delta-x_k\|<\varepsilon/3$ for any $n>n(\varepsilon,k)$. Since $k_n\rightarrow+\infty$ as $n\rightarrow+\infty$, for any $\varepsilon>0$, we have $\|x_{k_n}^\delta-x^*\|<\varepsilon$ as $k_n>k:=\max\{{k_1(\varepsilon),k_1(\varepsilon)}\}$, $n>n(\varepsilon,k)$.
\end{proof}

\providecommand{\href}[2]{#2}
\providecommand{\arxiv}[1]{\href{http://arxiv.org/abs/#1}{arXiv:#1}}
\providecommand{\url}[1]{\texttt{#1}}
\bibliographystyle{IEEEtran}
\bibliography{mybib}

\end{document}